\documentclass[11pt, reqno]{amsart}
 
 \usepackage{amsmath}
 \usepackage{amssymb}

\newcommand{\mysection}[1]{\section{#1}
      \setcounter{equation}{0}}

\newtheorem{theorem}{Theorem}[section]

\theoremstyle{definition}

\theoremstyle{remark}
\newtheorem{remark}{Remark}[section]

\newcommand\bbeta{\text{\raise-.2ex\hbox{$\bm{\beta}$}}}

\makeatletter
 \def\dashint{%  
 \operatorname%
 {\,\,\text{\bf--}\kern-.98em\DOTSI\intop\ilimits@\!\!}}
 \makeatother

\newcommand\bR{\mathbb{R}}

\begin{document}

\title[Quasi-convex functions]
{A note on quasi-convex functions}

\author{N.V. Krylov}
 
\email{nkrylov@umn.edu}
\address{127 Vincent Hall, University of Minnesota,
 Minneapolis, MN, 55455}



\begin{abstract}
We present an example of smooth quasi-convex
functions in the positive octant of $\bR^{3}$ which cannot be obtained as
the images of  convex smooth
functions under a monotone smooth mappings
of $\bR$.
 
\end{abstract}

\maketitle

\mysection{Introduction} 

Quasi-convex functions play an important role
in problems related to continuous optimization and mathematical
programming
such as generalizations of the von Neumann minimax theorem, the
Kuhn-Tucker saddle-point theorem, and other 
optimization problems related to consumer demand and 
indirect utility function (see, for instance, \cite{GP_71}, 
\cite{Mo_12}  and the references
therein).

According to \cite{GP_71} a real-valued function
$u(x)$ defined in a convex subset $E$ of
the Euclidean space $\bR^{d}$ is called quasi-convex
if 
$$
u(\lambda x+(1-\lambda)y)\leq \max[u(x),u(y)]
$$
as long as $\lambda\in[0,1]$, $x,y\in E$.
Convex functions in \cite{GP_71} are those
for which
$$
f(\lambda x+(1-\lambda)y)\leq \lambda f(x)+(1-\lambda)f(y).
$$

Numerous properties of quasi-convex functions
and their relation to convex functions
are discussed, for instance, in \cite{GP_71}, \cite{RV_73},
and \cite{BV-09}, in particular,
that $F[f(x)]$ is quasi-convex if $f$
is convex and $F$ is nondecreasing.
Somehow the following very
natural question is left untouched:
can any  quasi-convex function $u$
be represented as $F[f(x)]$ with convex $f$
and nondecreasing $F$? This issue is also avoided
in many other publications on the subject
of quasi-convexity. The reason for that is probably
because the answer to this question
is negative and the corresponding counterexamples
are given in \cite{DF_49}.
In these examples, however, $u$ is not smooth and for almost any point
of $E$ there exists a neighborhood
such that in that neighborhood the 
above representation still holds.

We want to present an example of a {\em smooth\/}
function $u$, which is quasi-convex
in a convex domain $E$ such that
there are  no smooth and strictly monotone functions $F$
such that $F[u(x)]$ is convex in a ball in $E$.
Then, of course, $u$ itself
 is not even locally an increasing smooth image
of a smooth convex function.
 This directly contradicts
the claim made in n.~1 of \cite{DF_49} that for any {\em smooth\/} 
 quasi-convex 
function $u$ one can find   strictly increasing $F$
such that  $F[u(x)]$ is convex.
Our arguments have much in common with Remark 5.14
of \cite{Kr_95}

Here is our result.

\begin{theorem}
                                        \label{theorem 10.14.1}
In $\bR^{3}$ consider the domain
$$
E:=\big\{(x,y,z):x,y,z>0 \big\},
$$
fix $\alpha>0$, and introduce the function
$$
u(x,y,z)=\frac{z^{\alpha}(x^{\alpha}+y^{\alpha})}
{x^{\alpha}y^{\alpha}}.
$$
Then $E$ is   convex, $u$ is quasi-convex in $E$,
and, if $\alpha\in(0,1]$, then
 for any $(x_{0},y_{0},z_{0})\in E$ and any twice continuously
differentiable function $F(t)$ on $\bR$ such that $F'[u(x_{0},y_{0},z_{0})]
\ne0$, the matrix of the second order
derivatives of $F[u( x,y,z)]$ at $(x_{0},y_{0},z_{0})$
is neither nonnegative nor
nonpositive. Hence,  
$F[u( x,y,z)]$ is neither convex nor concave
in any neighborhood of any point in $E$
if $F'>0$ or $F'<0$ on $\bR$.

\end{theorem}

Proof. First observe that the function
$$
v(x,y)=\frac{x^{\alpha}+y^{\alpha}}{x^{\alpha}y^{\alpha}}
=\frac{1}{x^{\alpha}}+\frac{1}{y^{\alpha}}
$$
is convex in $(0,\infty)^{2}$.
Then, in light of homogeneity of $u$    for $(x_{i},y_{i},z_{i})
\in(0,\infty)^{3}$, $\lambda_{i}\in[0,1]$, $i=1,2$,  
such that $\lambda_{1}+\lambda_{2}=1$, we have
$$
u\big(\lambda_{1}x_{1}+\lambda_{2}x_{2},
\lambda_{1}y_{1}+\lambda_{2}y_{2},
\lambda_{1}z_{1}+\lambda_{2}z_{2}\big)
=v\Big(\frac{\lambda_{1}x_{1}+\lambda_{2}x_{2}}
{\lambda_{1}z_{1}+\lambda_{2}z_{2}},
\frac{\lambda_{1}y_{1}+\lambda_{2}y_{2}}
{\lambda_{1}z_{1}+\lambda_{2}z_{2}}\Big) 
$$
$$
=v\Big(\mu_{1}\frac{ x_{1}}{z_{1}}+\mu_{2} \frac{ x_{2}}{z_{2}},
\mu_{1}\frac{ y_{1}}{z_{1}}+\mu_{2} \frac{ y_{2}}{z_{2}}\Big),
$$
where $\mu_{i}=\lambda_{i}z_{i}/(\lambda_{1}z_{1}+\lambda_{2}z_{2})$.
Since $\mu_{i}\geq0$ and $\mu_{1}+\mu_{2}=1$ and $v$ is convex,
the last expression above is less than
$$
\mu_{1}v\Big(\frac{ x_{1}}{z_{1}} ,
 \frac{ y_{1}}{z_{1}}\Big)+\mu_{2}v\Big(\frac{ x_{2}}{z_{12}} ,
 \frac{ y_{2}}{z_{2}}\Big)
=\mu_{1}u(x_{1},y_{1},z_{1})+
\mu_{2}u(x_{2},y_{2},z_{2})
$$
$$
\leq\max\big[u(x_{1},y_{1},z_{1}),u(x_{2},y_{2},z_{2})\big].
$$
This shows that $u$ is quasi-convex in $E$.

We now come to analyzing $F[u]$. Denote by $Dv$ the column-vector
gradient of $v$ and by $D^{2}v$ its matrix of the second-order
derivatives. By $a^{*}$ we mean the transpose of a matrix  $a$
and by $\langle a,b\rangle$ we mean the scalar product of
$a,b\in\bR^{3}$. We have
$$
D\{F[u]\}=F'[u]Du,\quad D^{2}\{F[u]\}= F''[u]Du(Du)^{*}+
F'[u]D^{2}u.
$$
We fix $(x_{0},y_{0},z_{0})\in E$ and take a column-vector
$\xi=(\xi_{1},\xi_{2},\xi_{3})\in\bR^{3}$ (written in a common abuse of notation
as a row vector) such that at $(x_{0},y_{0},z_{0})$ we have
$\langle Du,\xi\rangle=0$, which means that
$$
-\xi_{1}\alpha \frac{z^{\alpha}_{0}}
 {x_{0}^{\alpha+1} }-\xi_{2}\alpha \frac{z^{\alpha}_{0}}
 {y_{0}^{\alpha+1} }
+\xi_{3}\alpha\frac{z^{\alpha-1}_{0}( x^{\alpha}_{0}
+y^{\alpha}_{0} )}{x^{\alpha}_{0}y^{\alpha}_{0}}=0,
$$
\begin{equation}
                                                   \label{10.15.1}
  \xi_{3}=\frac{z_{0}x^{\alpha}_{0}y^{\alpha}_{0}}
{x^{\alpha}_{0}+y^{\alpha}_{0} }
\Big(\frac{\xi_{1}}{x^{\alpha+1}_{0}  }
 +\frac{\xi_{2}}{y^{\alpha+1}_{0}  }\Big).
\end{equation}

Furthermore,
$$
u_{xx}=\alpha(\alpha+1)
\frac{z^{\alpha}}{x^{\alpha+2}}
 ,\quad u_{xy}=0,\quad 
u_{xz}=-\alpha^{2}\frac{z ^{\alpha-1}}{x^{\alpha+1}},
\quad u_{yz}=-\alpha^{2}\frac{z ^{\alpha-1}}{y^{\alpha+1}},
$$
$$
u_{yy}=\alpha(\alpha+1)\frac{z^{\alpha}}{y^{\alpha+2}},
\quad u_{zz}=\alpha(\alpha-1)\frac{z^{\alpha-2} ( x^{\alpha} 
+y^{\alpha}  )}{x^{\alpha} y^{\alpha} },
$$
which implies that for $\xi$ satisfying \eqref{10.15.1} we have at
$(x_{0},y_{0},z_{0})$    %%%
$$
\big(F'[u]\big)^{-1}\langle \xi,D^{2}\{F[u]\}\xi\rangle
=\langle \xi,D^{2}u\xi\rangle
$$
$$
=\alpha(\alpha+1)
\frac{z_{0}^{\alpha}}{x_{0}^{\alpha+2}}\xi_{1}^{2}+
\alpha(\alpha+1)\frac{z_{0}^{\alpha}}{y_{0}^{\alpha+2}}\xi_{2}^{2}
$$
$$
-2\alpha^{2}\Big(\xi_{1}\frac{z_{0} ^{\alpha-1}}{x_{0}^{\alpha+1}}
+\xi_{2}\frac{z_{0} ^{\alpha-1}}{y_{0}^{\alpha+1}}\Big)
\frac{z_{0} x_{0}^{\alpha} y_{0}^{\alpha} }
{x_{0}^{\alpha} +y_{0}^{\alpha}  }
\Big(\frac{\xi_{1}}{x_{0}^{\alpha+1}   }
 +\frac{\xi_{2}}{y_{0}^{\alpha+1}  }\Big)
$$
$$
+\alpha(\alpha-1)\frac{z_{0}^{\alpha-2} ( x_{0}^{\alpha} 
+y_{0}^{\alpha}  )}{x_{0}^{\alpha} y_{0}^{\alpha} }
\frac{z_{0}^{2} x_{0}^{2\alpha} y_{0}^{2\alpha} }
{(x_{0}^{\alpha} +y_{0}^{\alpha} )^{2} }
\Big(\frac{\xi_{1}}{x_{0}^{\alpha+1}   }
 +\frac{\xi_{2}}{y_{0}^{\alpha+1}  }\Big)^{2}.
$$
Here the sum of the last two terms  equals
$$
z_{0}^{\alpha}\frac{  x_{0}^{\alpha} y_{0}^{\alpha} }
{x_{0}^{\alpha} +y_{0}^{\alpha}  }
\Big(\frac{\xi_{1}}{x_{0}^{\alpha+1}   }
 +\frac{\xi_{2}}{y_{0}^{\alpha+1}  }\Big)^{2}\big[\alpha(\alpha-1)-2\alpha^{2}
\big].
$$
Therefore, $\big(F'[u]\big)^{-1}\langle \xi,D^{2}\{F[u]\}\xi\rangle
=z_{0}^{\alpha}Q(\xi_{1},\xi_{2}) )$,
where
$$
Q(\xi_{1},\xi_{2} )=\frac{1}{x_{0}^{\alpha+2}}\Big[\alpha(\alpha+1)
-(\alpha^{2}+1) 
\frac{y_{0}^{\alpha}}{x_{0}^{\alpha} +y_{0}^{\alpha}  }\Big]\xi_{1}^{2}
$$
$$
+\frac{1}{y_{0}^{\alpha+2}}\Big[\alpha(\alpha+1)
-(\alpha^{2}+1) 
\frac{x_{0}^{\alpha}}{x_{0}^{\alpha} +y_{0}^{\alpha}  }\Big]\xi_{1}^{2}
$$
$$
-2(\alpha^{2}+1)\frac{  x_{0}^{\alpha} y_{0}^{\alpha} }
{x_{0}^{\alpha} +y_{0}^{\alpha}  }\frac{1}{x_{0}^{\alpha+1}}
\frac{1}{y_{0}^{\alpha+1}}\xi_{1}\xi_{2}
$$

It follows that to prove our claim it suffices to show that
the  quadratic form $Q$ with respect to $(\xi_{1},\xi_{2})$
is neither nonnegative nor nonpositive. For this to be true
we need to show that the determinant  of its matrix 
$$
\frac{1}{x_{0}^{\alpha+2}}
\frac{1}{y_{0}^{\alpha+2}}\Big[\alpha(\alpha+1)
-(\alpha^{2}+1) 
\frac{y_{0}^{\alpha}}{x_{0}^{\alpha} +y_{0}^{\alpha}  }\Big]
\Big[\alpha(\alpha+1)
-(\alpha^{2}+1) 
\frac{x_{0}^{\alpha}}{x_{0}^{\alpha} +y_{0}^{\alpha}  }\Big]
$$
$$
-(\alpha^{2}+1)^{2}\frac{  x_{0}^{2\alpha} y_{0}^{2\alpha} }
{(x_{0}^{\alpha} +y_{0}^{\alpha})^{2}  }\frac{1}{x_{0}^{2\alpha+2}}
\frac{1}{y_{0}^{2\alpha+2}}=:\frac{1}{x_{0}^{\alpha+2}}
\frac{1}{y_{0}^{\alpha+2}}R(\xi_{1},\xi_{2})
$$
is negative. One checks easily  that
$$
R(\xi_{1},\xi_{2})=\alpha^{2}(\alpha+1)^{2}-\alpha(\alpha+1)
(\alpha^{2}+1)= 
\alpha(\alpha+1)(\alpha-1),
$$
which is $<0$ if $\alpha<1$, and the theorem is proved in this case.

In case $\alpha=1$, take $\xi$ from above,
introduce $\eta=\xi+\kappa Du(x_{0},y_{0},z_{0})$,
and note that
$$
\big(F'[u]\big)^{-1}\langle \eta,D^{2}\{F[u]\}\eta\rangle
=\langle \xi,D^{2}u\xi\rangle+2\kappa\langle D^{2}u\xi,Du\rangle
$$
\begin{equation}
                                                        \label{10.17.1}
+\Big[F''[u]\big(F'[u]\big)^{-1}|Du|^{4}+\langle D^{2}u Du,Du\rangle\Big]
\kappa^{2},
\end{equation}
where at $(x_{0},y_{0},z_{0})$ as is easy to check
$$
\langle \xi,D^{2}u\xi\rangle=\frac{2z_{0}}{x_{0}+y_{0}}\Big(\frac{\xi_{1}}{x_{0}}
-\frac{\xi_{2}}{y_{0}}\Big)^{2}.
$$
This quantity vanishes if $\xi_{1}=tx_{0},\xi_{2}=ty_{0}$, $t\in\bR$,
and for the right-hand side of \eqref{10.17.1} not to change sign,
say for $\kappa=1$,
when $t$ runs through $\bR$ it is necessary to have
$\langle D^{2}uDu,\xi\rangle=0$ for those $\xi_{1},\xi_{2}$.

However, with this choice at $(x_{0},y_{0},z_{0})$ we have $\xi_{3}=tz_{0}$
and
$$
D^{2}u\xi=t\Big(  \frac{z_{0}}{x_{0}^{2}} ,
  \frac{z_{0}}{y_{0}^{2}} , -\frac{1}{x_{0} }-\frac{1}{y_{0} }\Big)
$$
$$
\langle D^{2}u\xi,Du\rangle=
t\Big(-\frac{z_{0}^{2}}{x_{0}^{4}}- \frac{z_{0}^{2}}{y_{0}^{4}}
-\big(\frac{1}{x_{0} }+\frac{1}{y_{0}}\big)^{2}\Big).
$$
This shows that for $\alpha=1$ as well, the right-hand side
of \eqref{10.17.1} has different signs for any fixed 
$(x_{0},y_{0},z_{0})$ if we vary $\xi$ and $\kappa$
and finishes the proof of the theorem. 

\begin{remark} 
One can show that,
if $\alpha>1$, then for any point $(x_{0},y_{0},z_{0})\in E$
one can find large $\lambda>0$ such that $\exp(\lambda u)$
is strictly convex in a neighborhood of $(x_{0},y_{0},z_{0})$.
Of course, $\lambda\to\infty$ as $\alpha\downarrow1$.
\end{remark}

\end{document}